\documentclass[12pt]{article}

\usepackage{amssymb,amsmath}
\usepackage{color}
\usepackage{latexsym,amssymb,amsfonts}
\usepackage{mathrsfs}
\usepackage{graphicx}
\usepackage{psfrag}
\usepackage{amsmath, amsbsy}
\usepackage{amsopn, amstext}
\usepackage{amsmath,multicol,amsopn}
\usepackage{latexsym,amssymb}

\def\sqr#1#2{{\vcenter{\vbox{\hrule height.#2pt
            \hbox{\vrule width.#2pt height#1pt \kern#1pt \vrule width.#2pt}
            \hrule height.#2pt}}}}
\def\signed #1{{\unskip\nobreak\hfil\penalty50
    \hskip2em\hbox{}\nobreak\hfil#1
    \parfillskip=0pt \finalhyphendemerits=0 \par}}
\def\endpf{\signed {$\sqr69$}}

\def\dbR{{\mathop{\rm l\negthinspace R}}}
\def\dbC{{\mathop{\rm l\negthinspace\negthinspace\negthinspace C}}}

\def\3n{\negthinspace \negthinspace \negthinspace }
\def\2n{\negthinspace \negthinspace }
\def\1n{\negthinspace }

\def\dbC{{\mathop{\rm l\negthinspace\negthinspace\negthinspace C}}}

\def\ds{\displaystyle}

\def\dbN{{\mathop{\rm l\negthinspace N}}}

\def\dbR{{\mathop{\rm l\negthinspace R}}}
\def\={\buildrel \triangle \over =}

\def\resp{{\it resp. }}
\def\a{\alpha}
\def\b{\beta}
\def\g{\gamma}
\def\d{\delta}
\def\e{\varepsilon}

\def\l{\lambda}

\def\n{\nabla}

\def\t{\times}
\def\f{\varphi}
\def\th{\theta}
\def\o{\omega}

\def\i{\infty}
\def\ns{\noalign{\ss} }

\def\G{\Gamma}
\def\D{\Delta}

\def\L{\Lambda}
\def\Si{\Sigma}

\def\O{\Omega}

\def\cC{{\cal C}}

\def\cF{{\cal F}}
\def\cG{{\cal G}}
\def\cH{{\cal H}}

\def\cJ{{\cal J}}
\def\cK{{\cal K}}

\def\cP{{\cal P}}

\def\cT{{\cal T}}

\def\cl{{\cal l}}
\def\no{\noindent}

\def\ms{\medskip}
\def\bs{\bigskip}
\def\q{\quad}
\def\qq{\qquad}
\def\hb{\hbox}

\def\min{\mathop{\rm min}}

\def\pa{\partial}

\def\cd{\cdot}

\def\div{\hbox{\rm div$\,$}}

\def\cl{\overline}

\def\Re{{\mathop{\rm Re}\,}}
\def\Im{{\mathop{\rm Im}\,}}

\def\|{\Big |}
\def\({\Big (}
\def\){\Big )}
\def\[{\Big[}
\def\]{\Big]}
\def\be{\begin{equation}}
\def\bel{\begin{equation}\label}
    \def\ee{\end{equation}}
\def\bt{\begin{theorem}}
    \def\bcd{\begin{condition}}
        \def\ecd{\end{condition}}
    \def\et{\end{theorem}}
\def\bc{\begin{corollary}}
    \def\ec{\end{corollary}}
\def\bde{\begin{definition}}
    \def\ede{\end{definition}}
\def\bl{\begin{lemma}}
    \def\el{\end{lemma}}
\def\bp{\begin{proposition}}
    \def\ep{\end{proposition}}
\def\br{\begin{remark}}
    \def\er{\end{remark}}
\def\ba{\begin{array}}
    \def\ea{\end{array}}
\def\ed{\end{document}}
\def\ns{\noalign{\ms}}
\def\ds{\displaystyle}

\def\square#1{\vbox{\hrule\hbox{\vrule height#1%
        \kern#1\vrule}\hrule}}
\def\rectangle#1#2{\vbox{\hrule\hbox{\vrule height#1%
        \kern#2\vrule}\hrule}}

\font\tenbb=msbm10 \font\sevenbb=msbm7 \font\fivebb=msbm5

\newfam\bbfam
\scriptscriptfont\bbfam=\fivebb \textfont\bbfam=\tenbb
\scriptfont\bbfam=\sevenbb

\def\oO{{\overline \O}}

\newtheorem{lemma}{Lemma}[section]
\newtheorem{remark}{Remark}[section]

\newtheorem{theorem}{Theorem}[section]
\newtheorem{corollary}{Corollary}[section]

\newtheorem{definition}{Definition}[section]
\newtheorem{proposition}{Proposition}[section]
\newtheorem{condition}{Condition}[section]

\makeatletter

\@addtoreset{equation}{section}
\makeatother

\begin{document}
%\begin{CJK*}{GBK}{kai}
\title{\bf Global Carleman Estimate  and State Observation Problem for Ginzburg-Landau Equation}

\author{Fangfang Dou\thanks{School of Mathematical Sciences, University of Electronic Science and Technology
of China, Chengdu, China. {\small\it
            E-mail:} {\small\tt fangfdou@uestc.edu.cn}.},  Zhonghua Liao\thanks{School
        of Mathematics, Sichuan University, Chengdu
        610064, China. {\small\it
            E-mail:} {\small\tt liaozhonghua@stu.scu.edu.cn}} \ and Xiaomin Zhu\thanks{School
        of Mathematics, Sichuan University, Chengdu
        610064, China. {\small\it
            E-mail:} {\small\tt zhuxiaomin@stu.scu.edu.cn}}}

\date{}

\maketitle

\begin{abstract}
In this paper, we prove a global Carleman estimate for  the complex
Ginzburg-Landau operator with a cubic nonlinear term in a bounded
domain of $\dbR^n, n=2,3$.  As applications, we study state
observation problems for the Ginzburg-Landau equation.
\end{abstract}
\bs

\no{\bf Key Words}. Ginzburg-Landau equation, Carleman estimate,
state observation problem, conditional stability. \bs

\no{\bf 2000 Mathematics Subject
Classification}.  93C20, 93B07.

\section{Introduction}\label{is1}

Let $T>0$, $\O\subset\dbR^n(n=2,3)$ be a given bounded domain with a
$C^2$ boundary $\Gamma$. Let $\o\subset \O$ be a nonempty open
subset and $\G_0\subset \G$  be a nonempty open subset. Put
$$Q\=(0,T)\t\O, \quad \Si\=(0,T)\t\Gamma, \q Q_\o\=(0,T)\t\o,\q \Si_0\=(0,T)\t\G_0.$$

Consider the following Ginzburg-Landau type equation
\begin{equation}\label{3ga}
y_t-(1+ib)\D y+(1+ic)|y|^2y=0 \qq\text {  in }  Q,\\
\end{equation}
where  $i=\sqrt{-1}$, $b, c\in\dbR$  characterize linear and
nonlinear dispersion, respectively.

Ginzburg-Landau equation is introduced to serve as a
phenomenological description  of superconductivity in \cite{GL1950}.
 Now it is used to described many physical phenomena,  such as
superconductivity, superfluidity, Bose-Einstein condensation to
liquid crystals and strings in field theory, and instability waves
(e.g., \cite{AL2001,RL2010}). It is one of the most-studied
nonlinear equations both in the physics and mathematics community.
Particularly, we refer the readers to \cite{GJL2020} for a detailed
introduction of the well-posedness of Ginzburg-Landau equations.

In this paper, we focus on establishing some Carleman estimates for
solutions to \eqref{3ga}.  Carleman estimates, which serve as
important tools in the study of unique continuation problems,
observability and controllability  problems,  and  inverse problems
for partial differential equations, have been investigated
extensively (see
\cite{FLZ,Hormander2009,IY,KL2021,LLR2022,Lerner2019} and the rich
references cited therein). Generally speaking, Carleman estimates
are established for linear partial differential operators. The usual
way to  hand semilinear equation is combining the Carleman estimate
for the linearized equation  and some Sobolev embedding theorem
(e.g., \cite{DZZ}). Nevertheless, such technique put some
restriction on the semilinear term, which is not satisfied by a
cubic term.

%The known literatures usually considered Carleman estimates for
%linear or semilinear operators, with Carleman weight functions
%satisfies ``pseudo--convex condition" in the sense of seminal work
%of L. H\"ormander \cite[Definition 8.6.1]{LH1}. For strong nonlinear
%operators, the common technique to the relevant problems is proving
%Carleman estimate for linearized operator and then applying the
%fixed-point argument. Under some restrictions, Carleman estimates
%for Ginzburg-Landau equations are studied. For instants, Fu
%\cite{Fu-jfa} proved a Carleman estimate for linearized
%Ginzburg-Landau operator (parabolic operator with complex
%coefficient), Rosier and Zhang \cite{RZ} obtained a Carleman
%estimate for semilinear Ginzburg Landau operator ``$ \pa_t -(a+ib)\D
%+d$"($ d \in L^\i(Q)$) which is closely related to the result in
%\cite{Fu-jfa} and yielded local null controllability for
%Ginzburg-Landau equation with strong nonlinear term ``$
%Ru+(1+i\beta)|u|^2u $" ($ R>0, ~\b $ are constants) by subtle
%fixed-point argument.

In this paper, we are endeavor to derive Carleman estimates for the
Ginzburg-Landau equation \eqref{3ga}. For convenience,  we define
the Ginzburg-Landau operator by
\begin{equation}\label{a1}
\cF y\= y_t-(1+ib)\D y+(1+ic)|y|^2y.
\end{equation}
We first establish a modified pointwise estimate for the
Ginzburg-Landau operator,  and then get the desired Carleman
estimates.

As an application of the Carleman estimates established in this
paper, we consider a state observation problem of the
Ginzburg-Landau equation, that is, can one determine the solution to
\eqref{3ga} (with suitable boundary condition) from some partial
measurement of the solution on $Q_0$ or $\Si_0$? More precisely, we
consider the following state observation problems:

\vspace{0.1cm}
\begin{itemize}
  \item The first one is the {\bf Identification Problem}: Is the solution
$y $  to \eqref{3ga} be determined uniquely by the observation $
y|_{Q_\o}$ (\resp $\ds\frac{\pa y}{\pa\nu}\Big|_{\Si_0}$)?

  \item
If the answer to the Identification problem is positive, then it is
natural to ask the  {\bf Conditional Stability Problem}: Assume that
two solutions $y$ and $\hat y$ (to the equation \eqref{3ga}) are
given. Let $y|_{Q_\o}$  and $\hat y|_{Q_\o}$  be the corresponding
observations. Can we find a positive constant $C$ such that
\begin{equation}\label{cs1}
|y-\hat y | \leq C(y,\hat y)|\!|y|_{Q_\o}- \hat y|_{Q_\o}|\!|, \;\;
\(\resp | y-\hat y | \leq C(y,\hat y)\|\!\| \frac{\pa
y}{\pa\nu}\|_{\Si_0}-\frac{\pa \hat y}{\pa\nu}\|_{\Si_0} \|\!\|\)
\end{equation}
with appropriate norms in both sides?
\end{itemize}
\begin{remark}
In the formulation of the Conditional Stability Problem, the
constant  $C(y,\hat y)$ in \eqref{cs1} means that we need some a
priori bound  for the solutions $y$ and $\hat y$. This is natural
since we deal with semilinear equation. However, the choice of the a
priori bound is subtle. Indeed, if one assume that
$|y|_{L^\infty(0,T;L^\infty(\O))}+|\hat
y|_{L^\infty(0,T;L^\infty(\O))}<M$ for some $M>0$, then the
derivation of \eqref{cs1} for the semilinear equation with cubic
term will be the same as the linear equation since
$$
|y^3| \leq M^2|y|,\q |\hat y^3| \leq M^2|\hat y|.
$$
In this paper, we do not assume such kind of condition.
\end{remark}
%

%\bel{3ga1}\left\{\ba{ll}\ds
%y_t-(1+ib)\D y+(1+ic)|y|^2y=0 &\mbox {  in $Q$ },\\
%\ns\ds\frac{\pa y}{\pa\nu}=h \mbox{ or } y=g &\mbox{ on $\Si$ },
%\ea\right.\ee
%with observable data $y|_{Q_\o}$,where $\omega$ be any non-empty open subset of $\Omega$, and $Q_\omega\triangleq (0,T)\times\omega$.

State observation problem in the above type is studied extensively
for linear PDEs (see \cite{BY2017,Fu-jfa,FLZ,Y} and the rich
references therein). However, the semilinear case with cubic term
attracts very little attention. To our best knowledge, there is no
published works addressing state observation problem for
\eqref{3ga}. The cubic nonlinear terms in the Ginzburg-Landau
equation will bring us many challenges in the proof that we cannot
simply imitate the results of linear equations to obtain the desired
conditional stability. As will be seen in Section \ref{ss5}, some
technical obstacles should be overcome.

The rest of this paper is organized as follows. In Section
\ref{sss2}, we state an internal and a boundary Carleman estimates
for the Ginzburg-Landau operators. These Carleman estimates will be
established in Section \ref{sec3}.   Finally, as applications of the
Carleman estimates given in Section \ref{sss2}, we solve the
conditional stability problem for the observability problem of
Ginzburg-Landau equation  in Section \ref{ss5}.

\section{Global Carleman estimate for Ginzburg-Landau operators}\label{sss2}

We first recall the following known results.
\bl\label{h31}  {\rm (\cite{IY}) } There is a real function
$\psi_1\in C^4(\oO)$ such that
$$\left\{\ba{ll}\ds\psi_1>0\ \mbox{ in }\O,\\
\ns\ds\psi_1=0 \ \mbox{ on } \G,\\
\ns\ds|\n\psi_1(x)|>0\ \mbox{ for all } x\in\oO\setminus{\o}.
\ea\right.$$
\el

Take a bounded domain $J\subset \dbR^n$  such that $\pa J \cap
\overline \Omega \subset \G_0$ and that $\widetilde \Omega = J\cup \Omega\cup\G$
has a $C^4$ boundary $\pa\widetilde \Omega$. Applying Lemma \ref{h31} to
the domain $\widetilde \Omega$ with choosing the set of the critical
points of $\psi_1$ belonging to $J$, we get the following result:
\bl\label{h31-0}   There is a real function
$\psi_2\in C^4(\oO)$ such that
$$\left\{\ba{ll}\ds\psi_2>0\ \mbox{ in }\O,\q |\n \psi_2(x)|>0 \mbox{ in } \O,\\
\ns\ds\psi_2=0 \q\mbox {on }  {\G\setminus\G_0},\q \frac{\pa\psi_2}{\pa\nu}\le 0,\ \ \forall x\in {\G\setminus\G_0}.
\ea\right.$$
\el

Now we introduce the  weight functions for the Carleman estimates.
Let $ \l>1 $ and $\mu
>1$. For $ j=1,2$, put
\bel{p33} \varphi_j(t,x)={e^{\mu\psi_j(x)}\over{t(T-t)}}, \q
\rho_j(t,x)={{e^{\mu\psi_j(x)}-e^{2\mu|\psi_j|_{C(\cl{\O
};\;\dbR)}}}\over{t(T-t)}},\q\ell_j=\l\rho_j,\q \th_j=e^{\ell_j}.
\ee

Denote
\bel{1027-2}\ba{ll}\ds
\a_1=-\frac{1}{1+b^2},\q \b_1=\frac{b}{1+b^2},\q \a_2=\frac{1+bc}{1+b^2},\q \b_2=\frac{c-b}{1+b^2},
\ea\ee
%and
%   \bel{1101-d}
% \g_1=\frac{1}{\a_1+i\b_1}=|\g_1|^{2}(\a_1-i\b_1),\q \g_2=\a_2+i\b_2.
% \ee
\bel{cp}\ba{ll}\nonumber\ds \cP y\=(\a_1+i \b_1)y_t+\D y \ea\ee
and
\bel{cg}\ba{ll}\ds \cG y\=\cP y-(\a_2+i \b_2)|y|^2y. \ea\ee
Then $\cF=-(1+ib)\cG$.   In the following, we consider the operator
$\cG$ instead of $\cF$ for the simplicity of notations.

We assume that
\begin{condition}\label{con1}
For some $ r_0\in (0,1) $ and $ \d_0 \in (0,1/8)$,
\bel{1102-f} |b|\le r_0<1,\qq\a_2>0,\qq |\b_2|\le \d_0\a_2. \ee
\end{condition}

\br\label{bc} From the first inequality of Condition \ref{con1} and
\eqref{1027-2}, we know $\b_1<\frac{r_0}{1+r_0^2}<\frac12$. Hence,
Condition \ref{con1} means that the imaginary part $b,c$ cannot be
large,  so do  $|b|$ and $|c|$. We believe that this is a technical
assumption. However, we do not know how to drop it. On the other
hand, this condition is satisfied in some important case
(e.g.,\cite{AL2001}).\er

We have the following internal and   boundary Carleman estimate  for
the  operator $\cG$.
\begin{theorem}\label{1.1}
Assume that Condition \ref{con1} holds. For all $y\in
C([0,T];L^2(\O))\cap L^2(0,T;H^1(\O))$ such that $\cG y\in
L^2(0,T;L^2(G))$ and $\ds{\pa y}/{\pa\nu}=0$ or $ \ds y= 0 $ on
$\Si$, there is a $ \mu_1>0 $ such that for all $ \mu\geq\mu_1 $,
one can find two constants $ C=C(\mu_1)>0 $ and $
\lambda_1=\lambda_1(\mu_1) $ so that for all $\l\ge\l_1$,  there
holds
\begin{equation}\label{5}
\ba{ll}
\ds\int_Q(\lambda\varphi_1)^{-1}\theta_1^2\big(|y_t|^2+|\Delta y|^2\big)dxdt
+\int_Q\big(\th_1^{2}|y|^6+\th_1^{2}|y|^2|\n y|^2\big)dxdt\\
\ns\ds +\l\mu^2 \int_Q\th_1^2\varphi_1\big(\lambda^2\mu^2\varphi_1^2|y|^2+|\n y|^2
+\l\varphi_1|y|^4\big)dxdt\\
\ns\ds\leq C\left[\int_Q\theta_1^2|\cG
y|^2dxdt+\lambda^2\mu^2\int_{Q_\o}\theta_1^2\varphi_1^2\big(\l\mu^2\varphi_1|y|^2+|y|^4\big)dxdt\right].
\ea
\end{equation}
\end{theorem}
Here and in the rest of this paper, we use $C$ to denote a generic
positive constant depending on $\O$, $T$, $\o$,  $b$  and $c$
(unless otherwise stated), which may change from line to line.

\begin{theorem}\label{1.1-0}
Assume that Condition \ref{con1} holds.  For all $y\in
C([0,T];L^2(\O))\cap L^2(0,T;H^1(\O))$ such that $\cG y\in
L^2(0,T;L^2(G))$ and $y=0$ on $\Si$, there is a $
\mu_2>0 $ such that for all $ \mu\geq\mu_2 $, one can find two
constants $ C^*=C^*(\mu_2)>0 $ and $ \lambda_2=\lambda_2(\mu_2) $
such that for all $\l\ge\l_2$, it holds that
\begin{equation}\label{5-0}\nonumber
\ba{ll}
&\ds\int_Q(\lambda\varphi_2)^{-1}\theta_2^2\big(|y_t|^2+|\Delta y|^2\big)dxdt
+\int_Q\big(\th_2^{2}|y|^6+\th_2^{2}|y|^2|\n y|^2\big)dxdt\\
\ns&\ds +\l\mu^2 \int_Q\th_2^2\varphi_2\big(\lambda^2\mu^2\varphi_2^2|y|^2+|\n y|^2
+\l\varphi_2|y|^4\big)dxdt\\
\ns&\ds\leq C^*\left(\int_Q\theta_2^2|\cG
y|^2dxdt\!+\l\mu\int_{\Si_0}\varphi_2\theta_2^2
\frac{\partial\psi_2}{\partial\nu}\|\frac{\partial
y}{\partial\nu}\|^2 d\G dt\right). \ea
\end{equation}
\end{theorem}

\br\label{r1} Recalling (\ref{1027-2}) for the definitions of $\a_2$
and $\b_2$, it is easy  to see that when $b\rightarrow 0$ and
$c\rightarrow0$, all the assumptions in \eqref{1102-f} are
satisfied. Thus one can get the Carleman estimate for the semilinear
heat operator ``$y_t-\D y+|y|^2y$" immediately. \er

\br\label{remark-0}  Carleman estimates for the heat operator
$\pa_t-\D$ and complex parabolic operator $\pa_t-(1+ib)\D$ with
Dirichlet boundary condition have been established in \cite{FI} and
\cite{Fu-jfa}, respectively. We mention that, following the proofs
of above two theorems in next section, we can easily obtain the
Carleman inequality for the parabolic operator with complex
coefficient $\pa_t-(1+ib)\D$ with Neumann boundary condition:
\begin{equation}\label{05-a}\nonumber
\ba{ll}
&\ds\int_Q(\lambda\varphi_1)^{-1}\theta_1^2\big(|y_t|^2+|\Delta y|^2\big)dxdt
+\l\mu^2 \int_Q\th_1^2\varphi_1\big(\lambda^2\mu^2\varphi_1^2|y|^2+|\n y|^2\big)dxdt\\
\ns&\ds\leq C\left[\int_Q\theta_1^2\big|y_t-(1+ib)\D
y\big|^2dxdt+\lambda^2\mu^2\int_{Q_\o}\varphi_1^3\theta_1^2\l\mu^2|y|^2
dxdt\right]. \ea
\end{equation}
and
\begin{equation}\label{05-b}\nonumber
\ba{ll}
\ds\int_Q(\lambda\varphi_2)^{-1}\theta_2^2\big(|y_t|^2+|\Delta
y|^2\big)dxdt
+\l\mu^2 \int_Q\th_2^2\varphi_2\big(\lambda^2\mu^2\varphi_2^2|y|^2+|\n y|^2\big)dxdt\\
\ns \ds\leq C^*\[\int_Q\theta_2^2\big|y_t-(1+ib)\D
y|^2dxdt+\l\mu\int_{\Si_0}\varphi_2\theta_2^2
\frac{\partial\psi_2}{\partial\nu}\|\frac{\partial
y}{\partial\nu}\|^2d\G dt\], \ea
\end{equation}
respectively.
\er

\section{Proof of Theorems \ref{1.1} and \ref{1.1-0}}\label{sec3}

In the following context, we denote by $\overline z, \Re z$ and $\Im
z$   the complex conjugate, real part and imaginary part of a
complex number $z\in\dbC$, respectively.

\subsection{A weighted identity for Ginzburg-Landau operator}\label{ss2}

In this subsection, we establish a pointwise weighted identity for
operator ``$\cG $" which is a key to prove our Carleman estimates.

\par Fix a weight function $\ell\in C^2(\dbR^{1+n};\dbR)$, and put
\bel{xz1802121}\nonumber
\th=e^{\ell},\quad
v=\th y.
\ee
Some elementary
calculations yield that
\bel{1027-5}\ba{ll}\ds
&\ds \th\cP y\=\th\big[(\a_1+i\b_1)y_t+\D y\big]\\
\ns  &\ds =(\a_1+i\b_1)(v_t-\ell_tv)+\D v-2\n\ell\cdot\n v+|\n\ell|^2v-\D\ell v\\
\ns &\ds=I_1+I_2,
\ea\ee
where
\bel{2pde3}\left\{\ba{ll}\ds
I_1\=i\b_1
v_t-\a_1\ell_tv+\D v+|\n\ell|^2v,\\
\ns\ds  I_2\=\a_1 v_t-i\b_1\ell_tv-2\n\ell\cdot\n v-\D\ell v.
\ea\right.\ee
Therefore,
\bel{1027-16}
\th\cG y=\th\cP y-\th(\a_2+i \b_2)|y|^2y=I_1+I_2-(\a_2+i\b_2)\th^{-2}|v|^2v=\cJ_1+\cJ_2,
\ee
where
\bel{j1j2}\left\{\ba{ll}\ds
\cJ_1=I_1-\frac{3\a_2}{4}|v|^2 v\theta^{-2},\\
\ns\ds
\cJ_2=I_2-\frac{\a_2}{4}|v|^2v\theta^{-2}-i \b_2|v|^2v\theta^{-2}.
\ea\right.\ee

We have the following weighted identity for  $\cG y$ defined by
(\ref{cg}).

    \begin{theorem}\label{3.1}
    Let $\a_j,\ \b_j\in\dbR$ ($j=1,2$).   Assume that $y,\;v\in C^2(\dbR^{1+n}; \;\dbC)$ and
    $\Psi,\;\Phi, \;\ell\in C^2(\dbR^{1+n};\dbR)$ satisfying
    \bel{1027-3}\nonumber\Psi+\Phi=-\Delta\ell.\ee
    Then, we have
\begin{eqnarray}\label{6}&&
2\Re(\theta\cG y\overline\cJ_1)+\(M+\frac{3}{8}\a_1\a_2\theta^{-2}|v|^4\)_t+\n\cdot\cH(v)\nonumber\\
&&=|\cJ_1|^2+|\cJ_1+\Phi v|^2+B|v|^2 +4\Re\sum_{j,k=1}^n\ell_{x_jx_k}v_{x_j}\overline{v}_{x_k}+2\Phi|\nabla v|^2+E\th^{-2}|v|^4+\frac{3\a_2^2}{8}\th^{-4}|v|^6 \nonumber\\
&&\q +\frac{\a_2}{4}\th^{-2}\big|\n
|v|^2\big|^2+\frac{\a_2}{2}\th^{-2}|v|^2|\n
v|^2+U+2\b_1\(\Phi+\frac{\a_2}{4}\th^{-2}|v|^2\)\Im(\overline{v}v_t),
\end{eqnarray}
    where $I_1$ and $I_2$ are given in (\ref{2pde3}), $\cJ_1$ and $\cJ_2$ are given in (\ref{j1j2}), in addition,
    \bel{1a3}\left\{\ba{ll}\ds
    M\=\left[(\a_1^2+\b_1^2)\ell_t-\a_1|\nabla\ell|^2\right]|v|^2+\a_1|\nabla v|^2-2\b_1\Im(\nabla\ell\cd\nabla\overline{v}v),\\
    \ns\ds V(v)\=4\Re (\nabla\ell\cd\nabla\overline{v})\nabla v-2|\nabla v|^2\nabla\ell-2\a_1\Re(\overline{v}_t\nabla v)+2\b_1\Im(\overline{v}_tv\nabla\ell)\\
    \ns\ds\qq\q +2\Im(\ell_t\overline{v}\nabla v)-2\Psi\Re(\overline{v}\nabla v)+2(|\nabla\ell|^2-2\a_1\ell_t)\nabla\ell|v|^2,\\
    \ns\ds
    B\=(\a_1^2+\b_1^2)\ell_{tt}+2\a_1\Phi\ell_t-4\a_1\nabla\ell\cd\nabla\ell_t+4\sum_{j,k=1}^n\ell_{x_jx_k}\ell_{x_j}\ell_{x_k}-2\Phi|\nabla\ell|^2-\Phi^2,
    \ea\right.\ee
    and
    \bel{key0}
    \left\{\ba{ll}\ds
    \cH(v)=V(v)-\frac{\a_2}{4}\theta^{-2}|v|^4\nabla\ell+\frac{\a_2}{2}\theta^{-2}|v|^2\Re(\overline{v}\nabla v),\\
    \ns\ds
    E=\frac{\a_2}{2}|\n\ell|^2+\a_2\D\ell-\frac{\a_1\a_2}{4}\ell_t+\frac{3\a_2}{2}\Phi,\\
    \ns\ds U=-4\b_1\Im(\nabla\ell_t\cd\nabla\overline{v}v)-2\Re(\nabla\Psi\cdot\nabla v \overline{v})-2\Re (i\b_2\overline {\cJ_1}\th^{-2}|v|^2v).
    \ea\right.
\end{equation}
\end{theorem}

\br In Theorem \ref{3.1}, $U$ can be regarded as low order terms
with respect to $v$  and $\n v$. It can be absorbed finally by the
energy terms   $|v|^2$ and $|\n v|^2$. The last term $\ds
2\b_1(\Phi+\a_2\tau\th^{-2}|v|^2)\Im(\overline{v}v_t)$  involves the
principal part $v_t$. Since it is related to the choice of $\Phi$,
we can deduce a modified estimate of this term in Section \ref{ss4}.
\er

\br As explained in \cite{FLZ}, in order  to keep more flexibility
in the sequel, we also introduce two auxiliary functions $\Psi\in
C^1(\dbR^{1+n};\dbR)$ and $\Phi\in C(\dbR^{1+n};\dbR)$. The choice
of them  will be given later. \er

We recall the following result, which is useful in the proof of
Theorem \ref{1.1}.

\bl\label{2l1} { \upshape(\cite[Theorem 1.1]{FLZ})} Under the
assumptions of Theorem \ref{2l1}, it holds that

\bel{2a2}\ba{ll}&\ds\ds
2\Re(\th\cP_1 y\overline
{I_1})+M_t+\n\cdot V(v)\\
\ns&\ds =|I_1|^2+|I_1+\Phi
v|^2+B|v|^2 +4\Re\sum_{j,k=1}^n\ell_{x_jx_k}v_{x_j}\overline{v}_{x_k}+2\Phi|\nabla v|^2\\
\ns&\ds\q  -2\Re(\nabla\Psi\cd\nabla v
\overline{v})-4\b_1\Im(\nabla\ell_t\cd\nabla\overline{v}v)+2\b_1\Phi\Im(\overline{v}v_t),
\ea\ee where $I_1$ is given in (\ref{2pde3}), $M,\  V$ and $B$  are
given in (\ref{1a3}). \el

\ms

%Mention that, the idea for prove lemma \eqr{2l1} be given in \cite[Theorem 2.1]{FLZ}is, by taking $\a=\a_1,\ \b=\b_1$ and $(a^{jk})_{n\t n}=I_n$ (the unit matrix), noting that $\Psi+\Phi=-\D\ell$, a short calculation yields the desired result immediately.\endpf

\noindent{\it Proof of Theorem \ref{3.1}.}
%{\it Step 1. }
By (\ref{1027-16}) and (\ref{j1j2}), we have
\begin{equation}\label{3.6}
\ba{ll}
&\ds2\Re(\theta\cG y\overline{\cJ_1})=2|\cJ_1|^2+2\Re (\overline {\cJ_1}\cJ_2)\\
\ns&\ds=2|\cJ_1|^2+2\Re (\overline {I_1}I_2)-\frac{3\a_2}{2}\Re (|v|^2 \overline v\theta^{-2}I_2) -\frac{\a_2}{2}\Re (\overline {I_1}\th^{-2}|v|^2v)\\
\ns&\ds\q -2\Re (i\b_2\overline {\cJ_1}\th^{-2}|v|^2v)+\frac{3\a_2^2}{8}\th^{-4}|v|^6.
\ea
\end{equation}
From (\ref{1027-5}),  (\ref{2pde3}) and (\ref{j1j2}), we find that
\bel{1027-4}\ba{ll}\ds 2\Re(\overline {I_1}I_2)= 2\Re(\theta\cP
y\overline {I_1})-2|I_1|^2, \ea\ee and \bel{1027-24}\ba{ll}\ds
&\ds |I_1+\Phi v|^2=|I_1|^2+2\Phi \Re (I_1\overline v)+\Phi^2|v|^2\\
\ns&\ds=|I_1|^2+2\Phi\Re(\cJ_1\overline v)+\Phi^2 |v|^2+\frac{3\a_2}{2}\Phi\theta^{-2}|v|^4\\
\ns&\ds=|I_1|^2+|\cJ_1+\Phi v|^2-|\cJ_1|^2+\frac{3\a_2}{2}\Phi\theta^{-2}|v|^4.
\ea\ee
Combining (\ref{2a2})--(\ref{1027-24}), we  obtain
\begin{equation}\label{1027-6}
\ba{ll}&\ds
 2\Re(\theta\cG y\overline\cJ_1)+M_t+\n\cdot V\\
\ns&\ds=|\cJ_1|^2+|\cJ_1+\Phi
v|^2+B|v|^2 +4\Re\sum_{j,k=1}^n\ell_{x_jx_k}v_{x_j}\overline{v}_{x_k}+2\Phi|\nabla v|^2\\
\ns&\ds\q  +\frac{3\a_2^2}{8}\th^{-4}|v|^6+\frac{3\a_2}{2}\Phi\theta^{-2}|v|^4-\frac{3\a_2}{2}\theta^{-2}|v|^2 \Re (\overline vI_2)\\
\ns&\ds\q -\frac{\a_2}{2}\th^{-2}|v|^2\Re (\overline {I_1}v)+U+2\b_1\Phi\Im(\overline{v}v_t).
\ea\ee

\ms

%{\it Step 2. }
Comparing (\ref{1027-6}) with (\ref{6}), we still need to deal with
$ ``\ds -\frac{3\a_2}{2}\theta^{-2}|v|^2\Re ( \overline v I_2)"$ and
$``\ds -\frac{\a_2}{2}\th^{-2}|v|^2\Re (\overline {I_1}  v)"$.
Clearly,
\bel{1027-7}\ba{ll}\ds
-2\Re (\theta^{-2}|v|^2 \overline vI_2)\\
\ns\ds=-2\Re \[\theta^{-2}|v|^2 \overline v(\a_1 v_t-i\b_1\ell_tv-2\n\ell\cdot\n v-\D\ell v)\]\\
\ns\ds=-2\a_1 \theta^{-2}|v|^2\Re ( \overline vv_t)+4\theta^{-2}|v|^2\Re ( \overline v\n\ell\cdot\n v)+2\D\ell\theta^{-2}|v|^4.
\ea\ee
By integrating by parts, the first two terms on the right hand side
of (\ref{1027-7}) are respectively
\bel{1027-9}\ba{ll}&\ds -2\a_1 \theta^{-2}|v|^2\Re ( \overline
vv_t)=-\a_1\th^{-2}|v|^2(|v|^2)_t=-\frac{1}{2}(\a_1\th^{-2}|v|^4)_t-\a_1\ell_t\th^{-2}|v|^4,
\ea\ee
and
\bel{1027-11}\ba{ll}\ds 4|v|^2\theta^{-2}\Re ( \overline
v\n\ell\cdot\n v)=\n\cdot (\th^{-2}|v|^4
\n\ell)+\th^{-2}(2|\n\ell|^2-\D\ell) |v|^4. \ea\ee
Substituting (\ref{1027-9}) and (\ref{1027-11})  into (\ref{1027-7})
yields
\bel{1027-12}\ba{ll}\ds
-\frac{3\a_2}{2}\Re (\theta^{-2}|v|^2 \overline vI_2)\\
\ns\ds=-\frac{3\a_1\a_2}{8}(\th^{-2}|v|^4)_t-\frac{3\a_1\a_2}{4}\th^{-2}\ell_t|v|^4+\frac{3\a_2}{4}\n\cdot (\th^{-2}|v|^4 \n\ell)\\
\ns\ds\q +\frac{3\a_2}{4}(2|\n\ell|^2-\D\ell)\th^{-2}
|v|^4+\frac{3\a_2}{2}\D\ell\theta^{-2}|v|^4. \ea\ee

Further,  from the definition of $I_1$ in (\ref{2pde3}) and by
integrating by parts, we obtain that
\bel{1027-13}\ba{ll}\ds
-\frac{\a_2}{2}\th^{-2}|v|^2\Re(\overline{I_1}v)\\
\ns\ds=-\frac{\a_2}{2}\th^{-2}|v|^2\Re\big[v(-i\b_1
\overline v_t-\a_1\ell_t\overline v+\D \overline v+|\n\ell|^2\overline v)\big]\\
\ns\ds=-\frac{\b_1\a_2}{2}\th^{-2}|v|^2\Im (v\overline
v_t)-\frac{\a_2}{2}(|\n\ell|^2-\a_1\ell_t)\th^{-2}|v|^4-\frac{\a_2}{2}\th^{-2}|v|^2\Re
(v\D \overline v). \ea\ee Similarly, there holds
\bel{1027-15}\ba{ll}\ds
-\frac{\a_2}{2}\th^{-2}|v|^2\Re (v\D \overline v)\\
\ns\ds=-\frac{\a_2}{2}\n\cdot\big[\th^{-2}|v|^2 \Re (v\n\overline v)\big]+\frac{\a_2}{2}\th^{-2}|v|^2|\n v|^2+\frac{\a_2}{4}\th^{-2}\big|\n |v|^2\big|^2\\
\ns\ds\q -\frac{\a_2}{4}\n\cdot\big(\th^{-2}\n\ell
|v|^4\big)-\frac{\a_2}{4}\th^{-2}(2|\n\ell|^2-\D\ell)|v|^4. \ea \ee

Finally, by substituting (\ref{1027-12})--(\ref{1027-15}) into
(\ref{1027-6}), we obtain \eqref{6} immediately.\endpf

\subsection{Proof of Theorems \ref{1.1} and  \ref{1.1-0}}\label{ss4}

\indent \par {\it Proof of Theorem \ref{1.1}.} The proof is long.
Thus we divide it into five steps.

\ms

{\it Step 1. } Put $ \ell=\ell_1$ and $\th=\th_1  $ in Theorem
\ref{3.1}. For $j,k=1,\cdots,n$, we have
\bel{33ss2}
\ell_{1t}=\l\rho_{1t},\q \ell_{1x_j}=\l\mu\varphi_1\psi_{1x_j},\q
\ell_{1x_jx_k}=\l\mu^2\varphi_1\psi_{1x_j}\psi_{1x_k}+\l\mu\varphi_1\psi_{1x_jx_k}.\ee
In addition, it is easy to get that
\bel{3a3ss2}
|\rho_{1t}|\le C e^{2\mu |\psi_1|_{C(\overline\O)}}\varphi_1^2,\q |\varphi_{1t}|\le
C\varphi_1^2,
\ee
and
\begin{equation}\label{b1}
\ba{ll}
\ds |\ell_{1tt}|\3n&\ds =\|\lambda\frac{e^{\mu\psi_1}-e^{2\mu|\psi_1|_{C(\bar{\O})}}}{t^2(T-t)^2}\[2\frac{(2t-T)^2}{t(T-t)}+2\]\|\\
\ns&\ds \leq C\lambda\varphi_1^2e^{2\mu|\psi_1|_{C(\bar{\O})}}(1+e^{-\mu\psi_1}\varphi_1)\\
\ns&\ds \leq C\lambda\varphi_1^3e^{3\mu|\psi_1|_{C(\bar{\O})}}. \ea
\end{equation}
Recalling that $\Phi+\Psi=-\D\ell_1$, by choosing $
\Psi=-2\l\mu^2\f_1|\n \psi_1|^2 $, we obtain
\bel{1027-21}
\Phi=\lambda\mu^2\varphi_1|\nabla\psi_1|^2-\l\mu\varphi_1\D\psi_1.
\ee
Next,  by (\ref{1027-2}), we see that
\bel{1102-a} -1<\a_1<0,\  |\b_1|\le \frac{1}{2}. \ee In the sequel,
for $ k\in\dbN $, we denote by $ O(\mu^k) $ a function of order $
\mu^k $ for large $ \mu $, and by $ O_\mu(\l^k) $ a function of
order $ \l^k $ for large $ \l $, with $\mu$ being fixed as a
parameter. Recalling (\ref{1a3}) for the definition of $B$, and by
using (\ref{33ss2})--(\ref{1102-a}), a short calculation shows that
\begin{equation}\label{18}\ba{ll}\ds  B\3n &\ds=(\a_1^2+\b_1^2)\ell_{1tt}+2\a_1\Phi\ell_{1t}-4\a_1\nabla\ell_1\cd\nabla\ell_{1t}
+4\sum_{j,k=1}^n\ell_{1x_jx_k}\ell_{1x_j}\ell_{1x_k}-2\Phi|\nabla\ell_1|^2-\Phi^2\\
\ns&\ds\geq
2\l^3\mu^4\varphi_1^3|\n\psi_1|^4+\lambda^3\varphi_1^3O(\mu^3)+\varphi_1^3O_\mu(\l^2),
\ea
\end{equation}
and that
\bel{1027-22}\ba{ll}\ds
E=\frac{\a_2}{2}|\n\ell_1|^2 +\a_2\D\ell_1-\frac{\a_1\a_2}{4}\ell_{1t}+\frac{3\a_2}{2}\Phi\ge \frac{\a_2}{2}\l^2\mu^2\varphi_1^2|\n\psi_1|^2-C\a_2\l\mu^2\varphi_1.
\ea\ee

Next, by choosing $ \Phi $ as in  (\ref{1027-21}), we have that
\bel{1027-23}\ba{ll}\ds
 4\Re\sum_{j,k=1}^n\ell_{1x_jx_k}v_{x_j}\overline{v}_{x_k}+2\Phi|\nabla v|^2\\
\ns\ds\ge 4\l\mu^2\varphi_1 |\n\psi_1\cdot\n v|^2+(2\lambda\mu^2\varphi_1|\nabla\psi_1|^2-C\l\mu\varphi_1)|\n v|^2.
\ea\ee
Further, noting that $\a_1<0$, we obtain that
\bel{1027-023}\ba{ll}\ds
U\3n&\ds=-4\b_1\Im(\nabla\ell_{1t}\cd\nabla\overline{v}v)-2\Re(\nabla v \overline{v})\cdot\nabla\Psi-2\Re (i\b_2\overline {\cJ_1}\th^{-2}|v|^2v)\\
\ns&\ds\ge-C\l^2\mu^4\varphi_1^3|\n\psi_1|^4|v|^2-C\mu^2\varphi_1|\n\psi_1|^2|\n v|^2-\frac{32\b_2^2}{\a_2^2}|\cJ_1|^2-\frac{1}{32}\a_2^2\th_1^{-4}|v|^6.
\ea\ee
Combining the pointwise identity (\ref{6}) given in Theorem \ref{3.1} with (\ref{1102-a})--(\ref{1027-023}), we conclude
\bel{1027-30}\ba{ll}\ds
 2|\theta_1\cG y|^2+\(M+\frac{3}{8}\a_1\theta_1^{-2}|v|^4\)_t+\n\cdot\cH\\
\ns\ds\ge 2\lambda^3\mu^4\varphi_1^3|\nabla\psi_1|^4|v|^2+4\l\mu^2\varphi_1 |\n\psi_1\cdot\n v|^2+2\lambda\mu^2\varphi_1|\nabla\psi_1|^2|\n v|^2+\frac{\a_2}{4}\th_1^{-2}\big|\n |v|^2\big|^2\\
\ns\ds\q +\frac{\a_2}{2}\th_1^{-2}|v|^2|\n v|^2+\frac{\a_2}{2}\l^2\mu^2\varphi_1^2|\n\psi_1|^2\th_1^{-2}|v|^4+\frac{11}{32}\a_2^2\th_1^{-4}|v|^6
 -\L_1\\
\ns\ds\q+2\b_1\(\Phi+\frac{\a_2}{4}\th_1^{-2}|v|^2\)\Im(\overline{v}v_t),
\ea\ee
where
\bel{1027-31}\ba{ll}\ds
\L_1=C(\lambda^2\mu^4\varphi_1^3-\lambda^3\varphi_1^3O(\mu^3)-\varphi_1^3O_\mu(\l^2))|v|^2\\
\ns\ds\qq \ +C(\l\mu\varphi_1+\mu^2\varphi_1)|\n v|^2+C\l\mu^2\varphi_1\th_1^{-2}|v|^4.
\ea\ee

\ms

{\it Step 2. }  We first estimate ``$\ds
2\b_1\(\Phi+\frac{\a_2}{4}\th_1^{-2}|v|^2\)\Im(\overline{v}v_t)$".

%\par We first estimate ``$2\b_1\Phi\Im(\overline{v}v_t)$".
For simplicity, put
\begin{equation}\label{1101-d}
    \g_1\=\frac{1}{\a_1+i\b_1}=|\g_1|^2(\a_1-i\b_1),\q\g_2\=\a_2+i\b_2.
\end{equation}
From (\ref{2pde3}) and (\ref{1027-16}), we have
\begin{equation}\label{49}
    v_t=\g_1\left( \theta_1\mathcal{G}y+\g_1^{-1}\ell_{1t} v-\Delta v+2\nabla\ell_1\cdot\nabla v-|\nabla\ell_1|^2v+\Delta\ell_1 v+\g_2|v|^2v\theta_1^{-2}\right).
\end{equation}
Noting that $ \Phi=\l\mu^2\f_1|\n \psi_1|^2-\l\mu\f_1\D \psi_1 $, we
get that
%\bel{310}\ba{ll}
%   2\b_1\Phi\Im(\overline{v}v_t)\\ \ns\ds=2\b_1\lambda\mu^2\varphi_1|\nabla\psi_1|^2\Im(\overline{v}v_t)-2\b_1\lambda\mu\varphi_1\Delta\psi_1\Im(\overline{v}v_t)\\
%   \ns\ds =2\beta_1\l\mu^2\varphi_1|\n\psi_1|^2\Im(\g\bar{v}\th\cG y)-2\b_1\l\mu^2\varphi_1|\n\psi_1|^2\Im(\g_1\bar{v}\D v)\\
%   \ns\ds \q+4\b_1\l\mu^2\varphi_1|\n\psi_1|^2\Im(\g_1\bar{v}\n\ell\cd\n v)-2\b_1\l\mu^2\varphi_1|\n\psi_1|^2\Im(\g)(|\n\ell|^2-\D\ell)|v|^2\\
%   \ns\ds\q+2\b_1\lambda\mu^2\varphi_1|\nabla\psi_1|^2\Im(\g_1\g_2)\th^{-2}|v|^4-2\b_1\lambda\mu\varphi_1\Delta\psi_1\Im(\overline{v}v_t)\\
%   \ns\ds\geq-C|\g_1|^{2}\lambda^2\mu^4\varphi_1^3|v|^2-|\theta\mathcal{G}y|^2 -2\b_1\lambda\mu^2\varphi_1|\nabla\psi_1|^2\Im\left(\g_1\overline{v}\Delta v\right)\\
%   \ns\ds\q +4\b_1\lambda^2\mu^3\varphi_1^2|\nabla\psi|^2\Im\left(\g_1\overline{v}\nabla\psi_1\cd\nabla v\right)+ 2\b_1^2|\g_1|^2\lambda^3\mu^4\varphi_1^3|\nabla\psi_1|^4|v|^2\\
%   \ns\ds\q+2\b_1\lambda\mu^2\varphi_1|\nabla\psi_1|^2\Im(\g_1\g_2)\th^{-2}|v|^4-(\lambda\mu\varphi_1)^{-1}|v_t|^2,
%   \ea\ee
\bel{310}\ba{ll} \ds
2\b_1\Phi\Im(\overline{v}v_t)=2\b_1\lambda\mu^2\varphi_1|\nabla\psi_1|^2\Im(\overline{v}v_t)
-2\b_1\lambda\mu\varphi_1\Delta\psi_1\Im(\overline{v}v_t)\\
\ns\ds\geq-C|\g_1|^{2}\lambda^2\mu^4\varphi_1^3|v|^2-|\theta_1\mathcal{G}y|^2 -2\b_1\lambda\mu^2\varphi_1|\nabla\psi_1|^2\Im\left(\g_1\overline{v}\Delta v\right)\\
\ns\ds\q  +4\b_1\lambda^2\mu^3\varphi_1^2|\nabla\psi_1|^2\Im\left(\g_1\overline{v}\nabla\psi_1\cd\nabla v\right)+ 2\b_1^2|\g_1|^2\lambda^3\mu^4\varphi_1^3|\nabla\psi_1|^4|v|^2\\
\ns\ds\q
+2\b_1\lambda\mu^2\varphi_1|\nabla\psi_1|^2\Im(\g_1\g_2)\th_1^{-2}|v|^4
-\beta_1^2\lambda^3\mu^3\varphi_1^3|\Delta\psi_1^4||v|^2-(\lambda\mu\varphi_1)^{-1}|v_t|^2.
\ea\ee

By some direct calculation, we have
\begin{equation}\label{key4}
\ba{ll}
-2\b_1\lambda\mu^2\varphi_1|\nabla\psi_1|^2\Im\left(\g_1\overline{v}\Delta v\right)\\
\ns\ds=-2\b_1\lambda\mu^2\n\cdot \left(\varphi_1|\nabla\psi_1|^2\Im(\g_1\overline{v}\n v)\right)+2\b_1\Im(\g_1) \lambda\mu^2\varphi_1|\nabla\psi_1|^2|\nabla v|^2\\
\ns\ds\q +2\b_1\lambda\mu^2\nabla(\varphi_1|\nabla\psi_1|^2)\cd\Im(\g_1\overline{v}\nabla v)\\
\ea
\end{equation}
and
\begin{equation}\label{key5}
\ba{ll}
4\b_1\lambda^2\mu^3\varphi_1^2|\nabla\psi_1|^2\Im\left(\g_1\overline{v}\nabla\psi_1\cd\nabla v\right)\\
\ns\ds=-2i\b_1|\g_1|^{2}\lambda^2\mu^3\varphi_1^2|\nabla\psi_1|^2\[(\a_1-i\b_1)\overline{v}\nabla\psi_1\cd\nabla v-(\a_1+i\b_1)v\nabla\psi_1\cd\nabla \overline v\]\\
\ns\ds=-2\b_1^2|\g_1|^2\lambda^2\mu^3\n\cdot \left( \varphi_1^2|\nabla\psi_1|^2|v|^2\nabla\psi_1\right)+2\b_1^2|\g_1|^2\l^2\mu^3\n\cdot \left( \varphi_1^2|\nabla\psi_1|^2\nabla\psi_1\right)|v|^2\\
\ns\ds\q +4\a_1\b_1|\g_1|^{2}\lambda^2\mu^3\varphi_1^2|\nabla\psi_1|^2\Im\left(\overline{v}\nabla\psi_1\cd\nabla v\right)\\
\ns\ds\geq -2\b_1^2|\g_1|^2\lambda^2\mu^3\n\cdot \left( \varphi_1^2|\nabla\psi_1|^2|v|^2\nabla\psi_1\right)+2\b_1^2|\g_1|^2\l^2\mu^3\n\cdot \left( \varphi_1^2|\nabla\psi_1|^2\nabla\psi_1\right)|v|^2\\
\ns\ds\q -4\alpha_1^2|\g_1|^2\l\mu^2\varphi_1|\nabla\psi_1\cdot\nabla
v|^2-\b_1^2|\g_1|^2\lambda^3\mu^4\varphi_1^3|\nabla\psi_1|^4|v|^2.
\ea
\end{equation}
%
%\par Next, we estimate$\ds \frac{1}{2}\b_1\a_2\th^{-2}|v|^2\Im(\overline{v}v_t)$.
By (\ref{49}), we know that
\bel{1101-e}\ba{ll}\ds
\frac{1}{2}\b_1\a_2\th_1^{-2}|v|^2\Im(\overline{v}v_t)\\
\ns\ds=
\frac{1}{2}\b_1\a_2\th_1^{-2}|v|^2\Im\big(\overline{v}\g_1 \theta_1\mathcal{G}y-\g_1\overline{v}\Delta v+2\g_1\nabla\ell_1\cdot\nabla v\overline{v}\big)\\
\ns\ds\q
-\frac{1}{2}\b_1\a_2\th_1^{-2}|v|^4\Im(\g_1)(|\nabla\ell_1|^2-\Delta\ell_1)
+\frac{1}{2}\b_1\a_2\th_1^{-4}|v|^6\Im(\g_1\g_2). \ea\ee
For the first term in the right hand side of \eqref{1101-e}, we have
\bel{1101-f}\ba{ll}\ds
\frac{1}{2}\b_1\a_2\th_1^{-2}|v|^2\Im(\overline{v}\g_1
\theta_1\mathcal{G}y)\ge
-\frac{1}{32}\a_2^2\th_1^{-4}|v|^6-2\b_1^2|\g_1|^2|\th_1\cG y|^2, \ea\ee
and
\bel{1101-g}\ba{ll}\ds
\frac{1}{2}\b_1\a_2\th_1^{-2}|v|^2\Im\big(-\g_1\overline{v}\Delta v+2\g_1
\nabla\ell_1\cdot\nabla v\overline{v}\big)\\
\ns\ds=-\frac{1}{2}\b_1\a_2\th_1^{-2}|v|^2\Im(\g_1\overline{v}\Delta v)+\b_1\a_2\th_1^{-2}|v|^2\Im(\g_1\nabla\ell_1\cdot\nabla v\overline{v})\\
\ns\ds=-\frac{1}{2}\b_1\a_2\div \big[\th_1^{-2}|v|^2\Im(\g_1\overline{v}\n v)\big]+\frac{1}{2}\b_1\a_2\Im(\g_1)\th_1^{-2}|v|^2|\n v|^2\\
\ns\ds\q +\b_1\a_2\th_1^{-2}\Re(\overline v\n v)\cdot\Im
(\g_1\overline v\n v). \ea\ee
Recalling the definition of $ \g_1 $
in  (\ref{1101-d}), it is easy to see that
\bel{1101-h}\ba{ll}\ds
\b_1\a_2\th_1^{-2}\Re(\overline v\n v)\cdot\Im (\g_1\overline v\n v)\\
\ns\ds=-\frac{1}{4}\b_1^2\a_2\th_1^{-2}|\g_1|^2\big|\n|v|^2\big|^2+\frac{1}{2}\b_1\a_2\a_1|\g_1|^2\th_1^{-2}\Im \big((\overline v\n v)^2\big)\\
\ns\ds\ge -\frac{1}{4}\b_1^2\a_2\th_1^{-2}|\g_1|^2\big|\n|v|^2\big|^2+\frac{1}{2}\b_1\a_2\a_1|\g_1|^2\th_1^{-2}|v|^2|\n v|^2.
\ea\ee
\par Now, substituting (\ref{310})--(\ref{1101-h}) into (\ref{1027-30}), and noting
that $\Im(\g_1)=-\b_1|\g_1|^2$, $\Im(\g_1\g_2)=(\a_1\b_2-\a_2\b_1)|\g_1|^2$ and
$\b_1^2|\g_1|^2<1$, we obtain that
\bel{1101-j}\ba{ll}\ds
(3+2\b_1^2|\g_1|^2)|\theta_1\cG y|^2+\(M+\frac{3}{8}\a_1\theta_1^{-2}|v|^4\)_t+\n\cdot\cK(v)\\
\ns\ds\ge (2-\b_1^2|\g_1|^2+ 2\b_1^2|\g_1|^2)\lambda^3\mu^4\varphi_1^3|\nabla\psi_1|^4|v|^2+4(1-\a_1^2|\g_1|^2)\l\mu^2\varphi_1 |\n\psi_1\cdot\n v|^2\\
\ns\ds \q+2(1-\b_1^2|\g_1|^2)\lambda\mu^2\varphi_1|\nabla\psi_1|^2|\n v|^2+\frac{1}{2}\a_2(1+\beta_1^2|\gamma_1|^2)\l^2\mu^2\varphi_1^2|\n\psi_1|^2\th_1^{-2}|v|^4+\cT\th_1^{-4}|v|^6\\
\ns\ds \q+\frac{1}{2}(1-\b_1^2|\g_1|^2+\b_1\a_1|\g_1|^2)\a_2\th_1^{-2}|v|^2|\n v|^2-\tilde{\L}_1-(\lambda\mu\varphi_1)^{-1}|v_t|^2,
\ea\ee
where
\begin{eqnarray}\ds
    &&\cK(v)=\cH(v)+2\b_1\lambda\mu^2\varphi_1|\nabla\psi_1|^2\Im(\g_1\overline{v}\n v)+2\b_1^2|\g_1|^2\lambda^2\mu^3 \varphi_1^2|\nabla\psi_1|^2|v|^2\nabla\psi_1\nonumber\\
    \ns\ds&&\qq \q\; +\frac{1}{2}\b_1\a_2\th_1^{-2}|v|^2\Im(\g_1\overline{v}\n v),\label{key9}\\
    \ns\ds&&\cT=\(\frac{5}{16}-\frac{1}{2}\b_1^2|\g_1|^2\)\a_2^2+\frac{1}{2}\b_1\a_2\a_1\b_2|\g_1|^2,\label{cf}
\end{eqnarray}
and
\begin{equation}\label{b2}\nonumber
    \ba{ll}
    \ns\ds \tilde{\L}_1=\L_1+C|\g_1|^2\l^2\mu^4\varphi_1^3|v|^2-2\beta_1\l\mu^2\varphi_1|\n\psi_1|^2\Im(\g_1\g_2)\th_1^{-2}|v|^4\\
    \ns\ds\qq \ -2\b_1\lambda\mu^2\nabla(\varphi_1|\nabla\psi_1|^2)\cd\Im(\g_1\overline{v}\nabla v)-2\beta_1^2|\g_1|^2\l^2\mu^3\n\cd(\varphi_1^2|\n\psi_1|^2\n\psi_1)|v|^2\\
    \ns\ds\qq \ +\frac{1}{2}\b_1^2|\gamma_1|^2\a_2(\l\mu^2\varphi_1|\n\psi_1|^2+\l\mu\varphi_1\D\psi_1)\theta_1^{-2}|v|^4+\beta_1^2\lambda^3\mu^3\varphi_1^3|\Delta\psi_1^4||v|^2.
    \ea
\end{equation}

Furthermore,  it is easy to check that
\bel{1102-b}\ba{ll}\ds
-\b_1^2|\g_1|^2+ 2\b_1^2|\g_1|^2>0,\q  1-\b_1^2|\g_1|^2= \frac{1}{1+b^2}=-\a_1,\\
\ns\ds  1-\a_1^2|\g_1|^2=\frac{b^2}{1+b^2}=-\a_1b^2,\q  1-\b_1^2\ge\frac{1}{1+b^2}=-\a_1,\\
\ns\ds 1-\b_1^2|\g_1|^2+\b_1\a_1|\g_1|^2=\frac{1-b}{1+b^2}\ge \frac{1-r_0}{(1+b^2)}=-(1-r_0)\a_1.
\ea\ee

By (\ref{1027-2}), (\ref{1102-f}), (\ref{1102-a}) and (\ref{cf}) and
noting that $|\g_1|^2=1+b^2<2$, we get that
\bel{1102-h}\ba{ll}\nonumber\ds
\cT\ds\ge  \(\frac{5}{16}-\frac{1}{8}|\g_1|^2\)\a_2^2-\frac{1}{4}\a_2|\b_2||\g_1|^2\ge\(\frac{5}{16}-\frac{1}{8}|\g_1|^2\)\a_2^2-\frac{1}{32}\a_2^2|\g_1|^2\\
\ns\ds\q=\frac{5}{32}\alpha_2^2(2-|\g_1|^2)>0. \ea\ee
By (\ref{1101-j}) and (\ref{1102-b}), we end up Step 2 with the
following inequality:
\bel{1102-d}\ba{ll}\ds
(3+2\b_1^2|\g_1|^2)|\theta_1\cG y|^2+\(M+\frac{3}{8}\a_1\theta_1^{-2}|v|^4\)_t+\n\cdot\cK(v)\\
\ns\ds\ge 2\lambda^3\mu^4\varphi_1^3|\nabla\psi_1|^4|v|^2\!-\!2\a_1\lambda\mu^2\varphi_1|\nabla\psi_1|^2|\n v|^2\!+\!\frac{\a_2}{2}(1\!+\!\b_1^2|\gamma_1|^2)\l^2\mu^2\varphi_1^2|\n\psi_1|^2\th_1^{-2}|v|^4\\
\ns\ds\q -\frac{1}{8}(1-b^2)\a_2^2\th_1^{-4}|v|^6
-\frac{1}{2}\a_1\a_2(1-r_0)\th_1^{-2}|v|^2|\n
v|^2-\tilde{\L}_1-(\lambda\mu\varphi_1)^{-1}|v_t|^2. \ea\ee

\ms

\par {\it Step 3.} Recalling  the choice of $\psi_1$ in Lemma \ref{h31}, we know that
$$ \min_{x\in {\O\setminus{\o_0}}}|\n \psi_1|>0.$$
For simplicity, we put
\bel{1102-0a}\nonumber
s_0\=\min_{x\in {\O\setminus
        {\o_0}}}\{|\n \psi_1|^2,\ |\n \psi_1|^4\}>0.
\ee
Noting that $\th_1(0,x)=\th_1(T,x)=0$,  $\a_1<0,\ \a_2>0$ and $
\b_1^2|\g_1^2|^2\leq 1$, by intergrating (\ref{1102-d}) over $ Q$,
we have
\bel{1102-0b}\ba{ll}\nonumber\ds
\l\mu^2 s_0\int_0^T\int_{\O\setminus{\o_0}}\varphi_1\(2\lambda^2\mu^2\varphi_1^2|v|^2-2\a_1|\n v|^2+\frac{1}{2}\a_2\l\varphi_1\th_1^{-2}|v|^4\)dxdt\\
\ns\ds -\frac{1}{8}\a_2\int_Q\((1-b^2)\a_2\th_1^{-4}|v|^6+4(1-r_0)\a_1\th_1^{-2}|v|^2|\n v|^2\)dxdt\\
\ns\ds\le 5\int_Q|\th_1\cG y|^2dxdt+\int_Q\n\cdot\cK(v) dxdt
+\int_Q\tilde{\L}_1dxdt+\int_Q(\lambda\mu\varphi_1)^{-1}|v_t|^2dxdt.
\ea\ee
Thus,
\bel{1102-0c}\ba{ll}\ds
\q\l\mu^2 s_0\int_Q\varphi_1\big(2\lambda^2\mu^2\varphi_1^2|v|^2-2\a_1|\n v|^2
+\frac{1}{2}\a_2\l\varphi_1\th_1^{-2}|v|^4\big)dxdt\\
\ns\ds\q -\frac{1}{8}\a_2\int_Q\((1-b^2)\a_2\th_1^{-4}|v|^6+4(1-r_0)\a_1\th_1^{-2}|v|^2|\n v|^2\)dxdt
\\
\ns\ds\le 5\int_Q|\th_1\cG y|^2dxdt+\int_Q\n\cdot\cK(v)dxdt+
\int_Q\tilde{\L}_1dxdt+\int_Q(\lambda\mu\varphi_1)^{-1}|v_t|^2dxdt\\
\ns\ds\q +\l\mu^2
s_0\int_{Q_{\o_0}}\varphi_1\big(2\lambda^2\mu^2\varphi_1^2|v|^2-2\a_1|\n
v|^2+\frac{1}{2}\a_2\l\varphi_1\th_1^{-2}|v|^4\big)dxdt. \ea\ee

\ms

{\it Step 4.} We estimate the boundary term $\int_Q\n\cdot\cK(v)
dxdt $ for the case that $y=0$ on $\Si$. Recalling that $ v=\th_1 y $,
we have $ v|_\Si= 0 $ and $ v_{x_j}=\frac{\pa v}{\pa \nu}\nu_j
(j=1,\cdots,n)$. Noting that $ \frac{\pa \psi}{\pa \nu} \le 0$ on $
\G $,  by (\ref{key9}),  we have
\begin{equation}\label{v1}
\ba{ll}
\ds \int_Q \n\cd \cK(v)dxdt\3n&\ds =\int_Q \n\cd V(v) dxdt\\
\ns&\ds = 2\l\mu \int_\Si \f_1 \frac{\pa \psi_1}{\pa \nu}\|\frac{\partial v}{\partial\nu}\|^2\big(\sum_{j,k=1}^n \nu_j\nu_k\big)^2 d\Gamma dt\le 0. \ea
\end{equation}
Noting that $\o_0\subset \subset \o$, we can choose a cut off
function $ \xi\in C_0^\infty (\omega;[0,1]) $ so that $ \xi\equiv1 $
on $ \omega_0$. Then,
\bel{1101-ds1}\ba{ll}\ds
\int_{Q}\xi^2\varphi_1 \th_1^2|\n y|^2dxdt\\
\ns\ds=-\Re\int_{Q}\xi^2\varphi_1 \th_1^2\overline y\D y dxdt
-\Re\int_Q\n(\xi^2\varphi_1\th_1^2)\cdot\n y\overline ydxdt\\
\ns\ds\le \e\int_Q\frac{1}{\l^2\mu^2\varphi_1}\theta_1|\D y|^2dxdt+\frac{C}{\e}\l^2\mu^2\int_{Q_\o}\varphi_1^3\th_1^2|y|^2dxdt+\frac{1}{2}\int_{Q}\xi^2\varphi_1 \th_1^2|\n y|^2dxdt\\
\ns\ds\q +C\l^2\mu^2\int_{Q_\o}\varphi_1^3\th_1^2|y|^2dxdt, \ea\ee
where $\e>0$ small enough.

Next, noting $\th_1(0,x)=\th_1(T,x)=0$, we have
\begin{eqnarray}\label{2131}
&&-2\Re\displaystyle\int_Q
\th_1^2\frac{1}{\lambda\varphi_1}(\a_1+i\b_1) y_t\cdot\D\overline
ydxdt \nonumber\\ &&=2\Re\displaystyle\int_Q\displaystyle
(\a_1+i\b_1)y_t\n \left(\th_1^2 \frac{1}{\lambda\varphi_1}\right)\cdot
\n \overline y dxdt +2\Re\displaystyle\int_Q\displaystyle\th_1^2
\frac{1}{\lambda\varphi_1}(\a_1+i\b_1)\n y_{t}\cdot \n\overline ydxdt \nonumber\\
&&=2\Re\displaystyle\int_Q\displaystyle (\a_1+i\b_1)y_t\n\left(\th_1^2
\frac{1}{\lambda\varphi_1}\right)\cdot \n\overline ydxdt
-\a_1\displaystyle\int_Q\displaystyle\left(\th_1^2
\frac{1}{\lambda\varphi_1}\right)_t |\n y|^2 dxdt\\
&&\leq\displaystyle\frac{1}{2}\displaystyle\int_Q
\th_1^2\displaystyle\frac{1}{\lambda\varphi_1}|(\a_1+i\b_1)y_t|^2dxdt+C\displaystyle\int_Q\th_1^2\lambda\mu^2\varphi_1|\nabla
y|^2 dxdt.\nonumber
\end{eqnarray}
From \eqref{2131}, we conclude that
\begin{equation}\label{381}
\ba{ll}\ds
\int_Q\frac{1}{\lambda\varphi_1}\theta_1^2\(|(\a_1+i\b_1)y_t|^2+|\Delta y|^2\) dxdt\\
\ns\ds=\int_Q\frac{1}{\lambda\varphi_1}\theta_1^2\|\cG
y+(\a_2+i\b_2)|y|^2 y\|^2 dxdt-2\Re\displaystyle\int_Q
\th_1^2\frac{1}{\lambda\varphi_1}(\a_1+i\b_1)
y_t\cdot\D\overline ydxdt\\
\ns\ds\leq\int_Q\frac{1}{\lambda\varphi_1}\left(|\th_1\cG
y|^2+|\a_2+i\b_2|^2\th_1^2|y|^6\right)dxdt+C\displaystyle\int_Q\th_1^2\lambda\mu^2\varphi_1|\nabla
y|^2 dxdt.\ea
\end{equation}
\par Combining (\ref{1102-a}) and  (\ref{1102-0c})--(\ref{381}), we find that for $ \e>0 $
sufficiently small,  there is a $\mu_0>0$ such that for all
$\mu\ge\mu_0$,  there exist two constants $\l_0=\l_0(\mu)$ and
$C>0$, so that for all $\l\ge\l_0$,  (\ref{5}) holds.

{\it Step 5.} In this step, borrowing some idea from
\cite{FI,Lu2013}, we deal with the case that $\ds \frac{\pa y}{\pa
\nu}= 0 $ on $ \Si $. Let
\bel{p33-0}\ba{ll}\nonumber\ds \tilde
\varphi_1(t,x)={e^{-\mu\psi_1(x)}\over{t(T-t)}},\q \tilde
\rho_1(t,x)={{e^{-\mu\psi_1(x)}-e^{2\mu|\psi_1|_{C(\cl{\O
};\;\dbR)}}}\over{t(T-t)}},\q \tilde \ell_1=\l\tilde \rho_1,\q\tilde
\th_1=e^{\tilde \ell_1}. \ea\ee
Set $\eta=\tilde\th_1 y$. Then, one can obtain  a similar inequality
as (\ref{1102-0c}) with $\th_1,\ \ell_1,\varphi_1$ and $\cK(v)$ replaced
by $\tilde\th_1,\ \tilde\ell_1,\ \tilde\varphi_1$ and $\cK(\eta)$,
respectively. Noting that $ 0 <\tilde \th_1\le \th,\  0<\tilde
\varphi_1<\varphi. $ Thus, we have
\bel{1102-0k}\ba{ll}\nonumber\ds
\hb{The left-hand side of (\ref{1102-0c})}\\
\ns\ds\le 10\int_Q|\th_1\cG y|^2dxdt+\int_Q\n\cdot\(\cK(v)+\cK(\eta)\) dxdt+2\int_Q\tilde{\L}_1dxdt\\
\ns\ds\q +\int_Q(\l\mu\tilde{\varphi_1})^{-1}|\eta_t|^2dxdt+\int_Q(\lambda\mu\varphi_1)^{-1}|v_t|^2dxdt\\
\ns\ds\q +2\l\mu^2
s_0\int_{Q_{\o_0}}\varphi_1\(2\lambda^2\mu^2\varphi_1^2|v|^2-2\a_1|\n
v|^2+\frac{1}{2}\a_2\l\varphi_1\th_1^{-2}|v|^4\)dxdt. \ea\ee

Noting that $\psi_1=0$ on the boundary $\Si$, and that $v=\th_1 y$
and $\eta=\tilde\th_1 y$, we see that
\bel{1102-0r}\left\{\ba{ll}\ds
\th_1=\tilde \th_1,\q \ell_1=\tilde\ell_1,\q \varphi_1=\tilde\varphi_1,\q v=\eta &\hb { on }\Si,\\
\ns\ds \ell_{1t}=\tilde\ell_{1t}, \q v_t=\eta_t,\q
\n\ell_1=-\n\tilde\ell_1,\q \n\psi_1=-\n\tilde\psi_1 &\hb{ on }\Si.
\ea\right.\ee
By (\ref{1102-0r}), noting that $\ds\frac{\pa y}{\pa\nu}=0$ on
$\Si$, $v=\th_1 y$ and $\eta=\tilde \th_1 y$,  it is easy to see
that
\bel{1102-0q}\ba{ll}\ds
\n v\cdot\nu+\n\eta\cdot\nu =\th_1 (\n y+\n\ell_1 y)\cdot\nu+\tilde\th_1 (\n y+\n\tilde\ell_1 y)\cdot\nu=0 &\hb { on }\Si,\\
\ns\ds \big(|\nabla v|^2\nabla\ell_1+|\nabla
\eta|^2\nabla\tilde\ell_1\big)\cdot\nu=2\th_1^2|\n\ell_1|^2\Re (\n
y\overline y)\cdot\nu=0 &\hb{ on }\Si, \ea\ee
and that
\bel{1102-0w}\ba{ll}\ds
\Re \big[(\nabla\ell_1\cd\nabla\overline{v})\nabla v+(\nabla\tilde\ell_1\cd\nabla\overline{\eta})\nabla \eta\big]\cdot\nu\\
\ns\ds=2\th_1^2\Re\big[(\n\ell_1\cdot\n \overline
y)\n\ell_1\cdot\nu+|\n\ell_1|^2(\n\overline y\cdot\nu)\big] y=0\q \hb
{on } \Si. \ea\ee
Recalling (\ref{1a3}), (\ref{key0}) and (\ref{key9}) for the
definitions of $V(v)$ ,$\cH$ and $ \cK $,  by
(\ref{1102-0r})--(\ref{1102-0w}), it is easy to see that
\bel{1102-0t}\nonumber \int_Q\n\cdot\(\cK(v)+\cK(\eta)\) dxdt=0. \ee
Further, recalling (\ref{1027-31}) for $\tilde{\L}_1$, by
(\ref{1102-a}), and noting that $v=\th_1 y$, we conclude that there
is a $\mu_0>0$ such that for all $\mu\ge\mu_0$, there exist two
constants $\l_0=\l_0(\mu)$ and $C>0$, so that for all $\l\ge\l_0$,
it holds that
\bel{1102-a0c}\ba{ll}\ds
\q\l\mu^2 \int_Q\tilde{\th}_1^2\tilde{\varphi_1}
\big(\lambda^2\mu^2\tilde{\varphi_1}^2|y|^2+|\n y|^2+\a_2\l\tilde{\varphi_1}|y|^4\big)dxdt\\
\ns\ds\q +\int_Q\big(\tilde{\th}_1^{2}|y|^6+\tilde{\th}_1^{2}|y|^2|\n y|^2\big)dxdt\\
\ns\ds\le C\[\int_Q|\th_1\cG y|^2dxdt+\int_Q\frac{1}{\lambda\mu\varphi_1}|y_t|^2dxdt+\int_Q \frac{1}{\l\mu\tilde{\varphi_1}}|y_t|^2dxdt\\
\ns\ds\q\q +\l\mu^2 \int_{Q_{\o_0}}\varphi_1\th_1^2
\big(\lambda^2\mu^2\varphi_1^2|y|^2+|\n
y|^2+\a_2\l\varphi_1|y|^4\big)dxdt\]. \ea\ee
On the other hand, we can obtain a similar result as (\ref{381})
with $ \varphi $ replaced by $ \tilde{\varphi_1} $ as follows:
\begin{equation}\label{b3}
\ba{ll}\ds
\int_Q\frac{1}{\lambda\tilde{\varphi_1}}\tilde{\theta_1}^2\big(|(\a_1+i\b_1)y_t|^2
+|\Delta y|^2\big) dxdt\\
\ns\ds=\int_Q\frac{1}{\lambda\tilde{\varphi_1}}\tilde{\theta_1}^2\big|\cG
y+(\a_2+i\b_2)|y|^2 y\big|^2 dxdt-2\Re\displaystyle\int_Q
\tilde{\th_1}^2\frac{1}{\lambda\tilde{\varphi_1}}(\a_1+i\b_1)
y_t\cdot\D\overline ydxdt\\
\ns\ds\leq\int_Q\frac{1}{\lambda\tilde{\varphi_1}}\big[|\tilde{\th_1}\cG
y|^2+|\alpha_2+\beta_2|^2\tilde{\th_1}^2|y|^6\big]dxdt
+C\displaystyle\int_Q\tilde{\th_1}^2\lambda\mu^2\tilde{\varphi_1}|\nabla
y|^2 dxdt\\
\ns\ds \leq C\int_Q\frac{1}{\l}\big[|\th_1\cG
y|^2+|\alpha_2+\beta_2|^2\th_1^2|y|^6\big] dxdt + C\int_Q\th_1^2\l\mu^2\varphi_1|\n
y|^2dxdt.\ea
\end{equation}
Finally, combining  (\ref{1102-a}), (\ref{1101-ds1}),  (\ref{381}),
(\ref{1102-a0c})  and (\ref{b3}), for $\e>0$ sufficiently small, we
can get the desired result immediately.\endpf \ms {\it Proof of
Theorem \ref{1.1-0}. } For $\th_2$ given by (\ref{p33}) with $\psi_2$
satisfying Lemma \ref{h31-0},   noting that $\psi_2=0$ on
$\G\setminus\G_0$, and by Lemma \ref{2l1}, proceeding exactly almost
the same argument as the proof of Theorem \ref{1.1}, one can get the
desired result.
\endpf

\section{Conditional stability}\label{ss5}

\indent\par In this section, we study the state observation for the
following Ginzburg-Landau equation \bel{3ga1}\left\{\ba{ll}\ds
y_t-(1+ib)\D y+(1+ic)|y|^2y=0 &\mbox {  in  } Q,\\
\ns\ds\frac{\pa y}{\pa\nu}=h \mbox{ or } y=g &\mbox{ on } \Si.
\ea\right.\ee

Let $\varepsilon\in (0,\frac{T}{2})$  and
\begin{equation}\label{vb1}
Q_\varepsilon\triangleq [\varepsilon,T-\varepsilon]\times\Omega.
\end{equation}

\bt\label{thm1} Assume that (\ref{1027-2}) and (\ref{1102-f}) are
fulfilled. There exists a constant $C>0$ such that for any solutions
$u_1,u_2\in C([0,T];L^2(\O))\cap L^2(0,T;H^1(\O))$ of (\ref{3ga1})
with $u_2\in L^\i ([0,T];L^6(\O))$,  it holds that
\begin{equation}\label{A211}
\int_{Q_\varepsilon}\big(|u_1-u_2|^2 + |\nabla u_1-\nabla
u_2|^2\big)dxdt \leq
\cC(u_2)\int_{Q_{\o}}(|u_1-u_2|^2+|u_1-u_2|^4)dxdt.
\end{equation}
where
$$
\cC(u_2)=C\Vert u_2\Vert_{L^\i ([0,T];L^6(\O))}^8.
$$
\et
\begin{remark}
As we explained in the introduction, Theorem \ref{thm1} provides a
conditional stability result, that is, if we know a priori that
$u_2\in L^\i ([0,T];L^6(\O))$, then we get an estimate for
$|u_1-u_2|$ on $Q_\e$ by $|u_1-u_2|$ on $Q_\o$. This is reasonable
since we deal with an equation with a cubic term. Further, the
condition that $u_2\in L^\i ([0,T];L^6(\O))$ is not very
restrictively. From the well-posedness result for \eqref{3ga1}, we
know this holds when the initial data of \eqref{3ga1} corresponding
to $u_2$ belong to $H_0^1(\O)$, then $u_2\in L^\i ([0,T];L^6(\O))$
(e.g., \cite[Chapter 3]{GJL2020}).
\end{remark}
\begin{remark}
From the proof of Theorem  \ref{thm1}, we see that the constant
$C\Vert u_2\Vert_{L^\i ([0,T];L^6(\O))}^8$ can be replaced by
$C\Vert u_1\Vert_{L^\i ([0,T];L^6(\O))}^8$, that is, we can assume
an a priori bound of $u_1$ to get the conditional stability.
\end{remark}
\begin{remark}
The constant $\e$ can be chosen to be any positive number.
Therefore, we know the answer to the Identification Problem is
positive. On the other hand, from the proof of Theorem  \ref{thm1},
we see that the constant $\cC(u_2)$ in \eqref{A211} tends to
infinity as $\e$ tends to $0$. This is reasonable since (\ref{3ga1})
is time irreversible.
\end{remark}
\begin{remark}
If we know a priori that the $L^\i ([0,T];L^6(\O))$-norm of the
state of the system (\ref{3ga1}) is smaller than a positive constant
$M$, then we see that the conditional stability result \eqref{A211}
can be presented as follows:
$$
\int_{Q_\varepsilon}\big(|u_1-u_2|^2 + |\nabla u_1-\nabla
u_2|^2\big)dxdt \leq CM^8\int_{Q_{\o}}(|u_1-u_2|^2+|u_1-u_2|^4)dxdt.
$$
\end{remark}

{\it Proof of Theorem  \ref{thm1}.} Let $u_1,u_2\in
C([0,T];L^2(\O))\cap L^2(0,T;H^1(\O))$, $u_2\in L^\i
([0,T];L^6(\O))$ and $u_1$, $u_2$ be the solutions
of (\ref{3ga1}).  Let $z\=u_1-u_2$. Then  we have%$$\partial_t (u_1\!-\!u_2)\!+\!(1\!+\!ib)\D (u_1\!-\!u_2)\!+\!(1\!+\!ic)|u_1\!-\!u_2|^2(u_1\!-\!u_2)\!=\!3(1\!+\!ic)(|u_1|^2u_2\!-\!u_1|u_2|^2)$$
%with the boundary condition
%$$\frac{\pa (u_1-u_2)}{\pa\nu}=0.$$
%   \begin{equation}\label{ga1_111}\nonumber
%       \left\{\ba{ll}\ds
%       \partial_t (u_1\!-\!u_2)\!-\!(1\!+\!ib)\D (u_1\!-\!u_2)\!+\!(1\!+\!ic)|u_1\!-\!u_2|^2(u_1\!-\!u_2)\!=\!3(1\!+\!ic)(|u_1|^2u_2\!-\!u_1|u_2|^2) \mbox{ in }Q,\\
%       \ns\ds\frac{\pa (u_1-u_2)}{\pa\nu}=0 \mbox{ or } u_1-u_2=0 \mbox{ on }\Si.
%       \ea\right.
%   \end{equation}
%For the sake of simplicity, we set $z=u_1-u_2$, then by fundamental computation we have
%
\begin{equation}\label{ga1_11}
\left\{\ba{ll}\ds
\partial_t z-(1+ib)\D z+(1+ic)|z|^2z=3(1+ic)(|z|^2u_2+z|u_2|^2) &\mbox {  in $Q$ },\\
\ns\ds\frac{\pa z}{\pa\nu}=0\mbox{ or } z=0 &\mbox{ on $\Si$ }.
\ea\right.
\end{equation}
Applying Theorem \ref{1.1} to \eqref{ga1_11}, we obtain that
\begin{equation}\label{517-1}
\ba{ll}
\ds \l\mu^2 \int_0^T\int_{\Omega}\th_1^2\varphi_1\(\lambda^2\mu^2\varphi_1^2|z|^2+|\n z|^2+\l\varphi_1|z|^4\)dxdt+\int_0^T\int_{\Omega}\(\th_1^{2}|z|^6+\th_1^{2}|z|^2|\n z|^2\)dxdt\\
\ns\ds\leq
C\left[\int_0^T\int_{\Omega}\theta_1^2(|z|^2u_2+z|u_2|^2)^2dxdt
+\lambda^2\mu^2\int_0^T\int_{\o}\varphi_1^2\theta_1^2(\l\mu^2\varphi_1|z|^2+|z|^4)dxdt\right].
\ea
\end{equation}
It follows from H\"older's inequality that
\begin{equation}\label{vg1}
\ba{ll}
\ds \int_0^T\int_{\Omega}\theta_1^2(|z|^2u_2+z|u_2|^2)^2dxdt\\
\ns\ds \le 2\int_0^T\int_{\Omega}\theta_1^2(|z|^4|u_2|^2+|z|^2|u_2|^4)dxdt\\
\ns\ds \le \epsilon_0\int_0^T\int_{\Omega}\theta_1^2 |z|^6dxdt
+C_{\epsilon_0}\int_0^T\int_{\Omega}\theta_1^2|z|^2|u_2|^4dxdt, \ea
\end{equation}
where $ \e_0 $ is a sufficiently small positive constant which can
be absorbed by the left  side of (\ref{517-1}).

By H\"older's inequality again, we have
\begin{equation}\label{517-2}
\ba{ll} \displaystyle{\int_0^T\int_{\Omega}\theta_1^2|z|^2|u_2|^4dxdt}
=\displaystyle{\int_0^T\int_{\Omega}|\theta_1 z|^2|u_2|^4dxdt}\\
\ns\ds\leq\int_0^T\left[\left(\int_{\Omega}|\theta_1 z|^6dx\right)^\frac{1}{3}\left(\int_{\Omega}|u_2|^6dx\right)^\frac{2}{3}\right]dt\\
\ns\ds\leq C\Vert u_2\Vert^4_{L^\i([0,T];L^6(\O))}\int_0^T\parallel
\theta_1 z \parallel^2_{H^1(\Omega)}dt. \ea
\end{equation}
Choosing $ \l>  C \Vert u_2\Vert^4_{L^\i([0,T];L^\infty(\O))}  $,
and gathering together (\ref{517-1})--(\ref {517-2}), we obtain that
\begin{equation}\label{517-4}\nonumber
\ba{ll}
\ds \l\mu^2 \int_0^T\int_{\Omega}\th_1^2\varphi_1\big(\lambda^2\mu^2\varphi_1^2|z|^2
+|\n z|^2+\l\varphi_1|z|^4\big)dxdt +\int_0^T\int_{\Omega}\big(\th_1^{2}|z|^6
+\th_1^{2}|z|^2|\n z|^2\big)dxdt\\
\ns\ds\leq
C\lambda^2\mu^2\int_0^T\int_{\o}\varphi_1^2\theta_1^2\big(\l\mu^2\varphi_1|z|^2+|z|^4\big)dxdt.
\ea
\end{equation}
\par On the other hand, noting that
\begin{equation}\label{517-6}
0<\th_1(\varepsilon,x)\leq \th_1(t,x)\leq \th_1\(\frac{T}{2},x\),\q\q
\forall (x,t)\in Q_\varepsilon,
\end{equation}
we deduce that
\begin{equation}\label{517-5}\nonumber
\ba{ll}
\ds\int_\varepsilon^{T-\varepsilon}\int_{\Omega}\th_1^{2}\big(|z|^2 + |\nabla z|^2 \big)dxdt\\
\ns\ds \leq  \lambda^3\mu^4\varphi_1^3\int_0^{T}\int_{\Omega}\th_1^{2}\big(|z|^2 + |\nabla z|^2 \big)dxdt \\
\ns\ds\leq
C\lambda^2\mu^2\int_0^T\int_{\o}\varphi_1^2\theta_1^2\big(\l\mu^2\varphi_1|z|^2+|z|^4\big)dxdt.
\ea
\end{equation}
This, together with \eqref{517-6}, implies \eqref{A211}.
\endpf

\br Similar to the proof of Theorem \ref{thm1}, one can give the
same analysis to obtain the conditional stability under boundary
observations, which is stated as follows.

Let $u_1, u_2\in C([0,T];L^2(\O))\cap L^2(0,T;H^1(\O))$ be solutions
of the following equation
\begin{equation}\label{710-1}
    \left\{\ba{ll}\ds
    y_t-(1+ib)\D y+(1+ic)|y|^2y=0  &\mbox { in }Q,\\
    \ns\ds y=g  &\mbox{ on  }\Si,
    \ea\right.
\end{equation}
and  $u_2\in L^\infty(0,T;L^6(\O))$,  then it holds that
\begin{equation}\label{710-2}\nonumber
\ba{ll} \displaystyle{\int_{Q_\varepsilon}\big(|u_1-u_2|^2 + |\nabla
u_1-\nabla u_2|^2\big)dxdt}  \leq C\displaystyle{\int_{\Si_0}
\|\frac{\pa u_1}{\pa \nu}-\frac{\pa u_2}{\pa \nu}\|^2 d\G dt}. \ea
\end{equation}
\er

\section*{Acknowledgement}

This work is partially supported by the NSF of China (No. 11971333, 11931011, 12071061,
11971093),  the Science Fund for Distinguished Young Scholars of Sichuan Province (No. 2022JDJQ0035), the Applied Fundamental Research Program of Sichuan Province (No. 2020YJ0264), the Science Development Project of Sichuan University (No. 2020SCUNL201), and the Fundamental Research Funds for the Central Universities of UESTC (No. ZYGX2019 J094).

%\end{CJK*}
\end{document}